\nonstopmode 
\magnification=\magstep 1 
\input amstex

\documentstyle{amsppt}
\loadbold \loadmsbm\loadeufm 
\topmatter
\title ORDINARY SUBVARIETIES OF CODIMENSION ONE\endtitle
\author Ferruccio Orecchia\endauthor
\address Ferruccio Orecchia, Dipartimento di Matematica e Applicazioni "R. Caccioppoli", Complesso 
Universitario di Monte S. Angelo
- Via Cintia , 80126 Napoli-Italy \endaddress
\email ORECCHIA\@MATNA2.DMA.UNINA.IT.\endemail\
\thanks Work partially supported by 
MURST\endthanks 

\abstract In this paper we extend the properties 
of ordinary  points of curves [10] to ordinary closed points of one-dimensional affine 
reduced schemes and then to  ordinary  subvarieties of codimension one. 
\endabstract
\endtopmatter

\document
\head Introduction.\endhead 
 Let $A$ be the local ring, at a  point $x$,
 of an algebraic reduced variety $X$ over an algebraically closed field $k$.
Denote with 
 $e(A)=e$ the multiplicity of $X$ at $x$. Let  $\frak m$ be the  maximal ideal of $A$.
$Spec(G(A))$ is the tangent cone and $Proj(G(A))$ the projectivized tangent cone
to $X$ at $x$. Furthermore, if $G(A)_{red}=Spec(G(A)/nil(G(A))$, then $Spec(G(A)_{red})$
is the tangent cone to $X$ at $x$, viewed as a set. If $\overline A$ is the
normalization
of $A$, the branches of  $X$ at $x$
 are the points of 
$Spec(\overline A/\frak m\overline A)$.\medskip
 If $A$ is the local ring, at a point $x$,  of a curve $C$ and $\frak p_i$, 
$i=1,...,n$, are the minimal primes
of $G(A)$ then $Spec(G(A)/\frak p_i)$ are the tangents of $C$ at $x$ (that is the 
tangents are the lines of $Spec(G(A)_{red})$ and correspond to the points of  $Proj(G(A))$). 
 It is well known [10, Lemma-Definition 2.1] that each branch
 has a tangent
and that the following conditions are
equivalent:\medskip
 (a) The scheme $Proj(G(A))$ is reduced;\medskip
 (b) $Proj(G(A))$ consists of $e$ points;\medskip
 (c) there are $e=e(A)$ tangents at $x$;\medskip 
 (d) the curve $C$ has, at $x$, $e$ linear branches with distinct tangents.\medskip
A point $x$ which satisfies these
 conditions is said to be {\sl an ordinary point} . Any nonsingular point 
of the curve ($e=1$) is a trivial example of ordinary point. Generally,
 but not always [9, Section  4],
 an ordinary point has reduced tangent cone, in the sense that, if its 
tangents are in generic position,  then the tangent cone to $C$ at $x$ (i.e. $G(A)$) 
is reduced  [10, Theorem 3.3]. 
\medskip
In this paper we extend these results to subvarieties of codimension one.\medskip 
For this  we need first to extend the notion of ordinary point of a curve
 to the notion
of ordinary  closed point of a one-dimensional reduced local affine scheme $Spec(A)$.
In fact let  $\frak m$  be the maximal ideal of $A$ and $K$ be the algebraic closure of the 
residue  field $k(\frak m)$ of $A$.
Then the closed point of $Spec(A)$ is said to be {\sl  ordinary } if
 $Proj(G(A) \otimes_{k(\frak m)}  K)$ is reduced [Definition 2.4].

 We prove that the previous equivalences $(a)\Leftrightarrow(b)\Leftrightarrow(c)\Leftrightarrow(d)$
 continue to hold in this more
 general setting if we suitably define (geometric) tangents and (geometric) branches and we prove 
that if the ordinary closed
point of $Spec(A)$ has (geometric) tangents in generic position then $G(A) \otimes_{k(\frak m)}  K$ is reduced.\medskip

 Then, basing on this, using  results of [2] and the notion of normal flatness,
 we extend  the properties
 of ordinary points  to subvarieties of
codimension one.
Let $X=Spec(R)$ be an equidimensional algebraic variety and $Y=Spec(R/\frak q)$ be
 an irreducible codimension one subvariety  of $X$ of multiplicity $e=e(R_\frak q)$.
 If $A$ is the local ring of $X$ at a closed point $x$ of $Y$ then $x$ 
 is said to be an {\sl ordinary point}  if $Spec(G(A)_{red})$ 
 is the union of $e$ distinct linear spaces.  
$Y$ is said to be an {\sl ordinary subvariety} of $X$ if  the closed point of 
the one-dimensional scheme $Spec(R_\frak q)$ is ordinary. 
Assume that $Y$ is nonsingular at $x$ and    
$X$ is normally flat along $Y$ at $x$, then $x$ is an ordinary point of $X$ 
if and only if $X$ has, at $x$, $e$ linear branches with distinct tangent spaces.
If these tangent spaces are in generic position then $G(A)$ is reduced [Theorem 3.4]. 
Moreover $Y$ is an  ordinary subvariety of $X$ if and only if there exists an
open non-empty subset $U$ of $Y$ such that every closed point $x$ of $U$ is ordinary 
[Theorem 3.5].\medskip
This paper has been motivated by  [12] in which
 the conductor of a
variety with  ordinary (multiple) subvarieties of codimension
one is computed. We would like to thank  G. Castaldo for pointing out the paper
 [2].\medskip
Throughout the paper all rings are supposed to be commutative, with identity and 
noetherian.\medskip
Let $A$ be a local ring with maximal ideal $\frak m$.
With $H^0(A,n)=dim_k(\frak m^n/ \frak m^{n+1})$, $n\in \Bbb N$, we denote the
 Hilbert function of $A$ and with $e(A)$ the multiplicity of $A$ at $\frak m$. 
The embedding dimension $emdim(A)$ of $A$ is given by $H^0(A,1)$.
With $A^h$ we denote the henselization of $A$ and with $A^{hs}$ a fixed strict 
henselization of $A$ (see [6, Number 32, Section 18] or [13]).
 If $\frak p$ is an  ideal of $A$, 
$G_\frak p (A)= \bigoplus_{n\geq 0} (\frak p^n/\frak p^{n+1})$ 
is the associated 
graded ring with respect to $\frak p$. By $G(A)$ 
we denote the associated graded ring 
with respect to  $\frak m$. $Spec(G(A)$ is the {\sl tangent cone} to $Spec(A)$ at its closed point.  
If $S=\bigoplus_{n\geq 0} S_n$ is a standard graded finitely generated
algebra over a field $k$ ,
with maximal homogeneous ideal $\frak n$,  
$H(S,n)=dim_k S_n=H(S_\frak n,n)$
 denotes the Hilbert function of $S$ and  $emdim(S)=H(S,1)=emdim(S_\frak n)$
denotes the embedding dimension of $S$.  The multiplicity of
$S$ is  $e(S)=e(S_\frak n)$ . On has $e(A)=e(G(A))$ and $emdim(A)=emdim(G(A))$.
If $k(\frak m)$ is the residue field of $A$, 
$Hom_{k(\frak m)}(\frak m/\frak m^2, k(\frak m))$ is the {\sl tangent space}
of $Spec(A)$ at  its closed point. If $A$ is regular then the tangent cone and the tangent space 
of $A$ coincide.
\medskip If $R$ is any ring $dim(R)$ denotes the dimension of $R$.\medskip
\head 1. Generalities on branches and normal flatness\endhead
\definition {Definition 1.1} Let $A$ be a local reduced ring
 with maximal ideal $\frak m$. Let $K$ be the algebraic closure of the residue field $k(\frak m)=A/\frak m$ and
$\overline A$ be the normalization of $A$. A {\sl branch}   of
 $Spec(A)$ at its closed point is a  point of $Spec(\overline A/\frak m\overline A)$. A {\sl 
geometric branch}  of $Spec(A)$ at its closed point is a  point of 
$Spec(\overline A/\frak m\overline A\otimes_{k(\frak m)} K)$. If $X$ is a scheme a
{\sl branch} (respectively a {\sl geometric branch})  of $X$ at a closed point $x$
 is a branch (respectively a geometric branch) at the closed point of $Spec(A)$,
where $A$ is the local ring of $X$ at $x$.\enddefinition
\remark{Remark} Since $\overline A$ is integral over $A$ the ring $\overline A/\frak m\overline A$
is zero-dimensional and then the ring 
$\overline A/\frak m\overline A\otimes_{k(\frak m)} K$ 
is zero-dimensional [6, $N^0$ 24, Corollaire 4.1.4]. Then the points of the schemes
$Spec(\overline A/\frak m\overline A)$ and 
$Spec(\overline A/\frak m\overline A\otimes_{k(\frak m)} K)$ are closed.\endremark
\proclaim{Lemma 1.2} If $A$ is a local reduced ring, there is a canonical bijection 
between the set 
of branches of $Spec(A)$ and the set of minimal primes of $A^h$ and there is a canonical 
bijection between the set of geometric branches of $A$ and the set of minimal primes
of $A^{sh}$.\endproclaim
\demo{Proof} [5, Theorem 2.1]\enddemo
\definition{Definition 1.3} Let $A$ be a local reduced ring and $\gamma$ 
 be a branch (respectively a geometric branch) of $Spec(A)$
at  its closed point.
 If $\frak p$ is the minimal prime of $A^h$ (respectively $A^{hs}$)
 corresponding to $\gamma$,  the {\sl order of the branch} (respectively
 the {\sl order of the geometric branch}) $e(\gamma)$ is the multiplicity 
$e(A^h/\frak p)$ (respectively $e(A^{sh}/\frak p)$). The branch (respectively 
the geometric branch) $\gamma$ is {\sl linear} if $A^h/\frak p$
 (respectively $A^{hs}/\frak p$)
is a regular ring.\enddefinition
\remark{Remark} The order of a geometric branch is independent from the strict henselization
of $A$ chosen [1, Corollary 1.9]\endremark
\proclaim{Proposition 1.4} Let $A$ be a reduced local ring of dimension one or the
local ring at a closed point of an algebraic equidimensional reduced variety (over
an algebraically closed field $k$).
Let $\gamma_i$, $i=1,...,n$ be the branches (respectively the geometric branches) of
$Spec(A)$ at the closed point. Then:\medskip
 (a) $e(A)=\sum_{i=1}^ne(\gamma_i)$;\medskip
 (b) $\gamma_i$ is linear if and only if its order is one;\medskip
 (c) $n=e(A)$ if and only if all the branches are linear.\endproclaim
\demo{Proof} $(a)$, $(b)$ (see  [1, Propositions 2.7 and 2.10]). 
$(c)$ is obvious consequence of $(a)$ and $(b)$.\qed\enddemo

\definition{Definition 1.5} Le $\frak p$ be a prime ideal of a local ring $A$.
$A$ is {\sl normally flat along $\frak p$ } if $\frak p^n/\frak p^{n+1}$
is flat over $A/\frak p$, for all $n\geq 0$. Note that $A$ is normally flat
along $\frak p$ if and only if $G_\frak p(A)$ is free over $A/\frak p$.\enddefinition

\proclaim{Theorem 1.6} Let $A$ be a local ring and $\frak p$ a prime ideal of $A$
such that $A/\frak p$ is regular of dimension $d$ and $A$ is normally flat along $\frak p$.
Then:\medskip
 (a) there is an isomorphism of graded $k$-algebras
  
 $$G(A) \cong (G_\frak p (A) \otimes _{A/ \frak p} k)[T_1,...,T_d]$$\medskip

(b) the following rings have the same Hilbert functions:\medskip 
  $G_\frak p (A) \otimes _{A/ \frak p} k$ (over $k$) and $G(A_\frak p)$
(over $k(\frak p))$\medskip.
\endproclaim
\demo{Proof} $(a)$ [7, Corollary (21.11)].\medskip

$(b)$ [8, Chapter II, Corollary 2].\qed\enddemo

Let $X=Spec(R)$ be a reduced variety over an algebraically closed
field $k$ and $Y=Spec(R/\frak q)$ be an irreducible  subvariety $Y$ 
of $X$.

\definition {Definition 1.7}If $x$ is any closed point of $Y$, $A$ 
is the local ring of $X$ 
at $x$ and $\frak p=\frak qA$ is 
the prime ideal in
 $A$ defining the subvariety $Y$, then:\medskip

$(i)$ $X$ is {\sl normally flat along $Y$ at $x$} if 
$A$ is normally flat along $\frak p$;
\medskip

$(ii)$ $Y$ is 
{\sl nonsingular
at} $x$ if $A/\frak p$ is regular.\enddefinition

\proclaim{Theorem 1.8} There exists an open nonempty subset $U$ of $Y$ such that,
 for every closed point $x$
of $U$,  $Y$ is nonsingular at $x$ and
 $X$ is normally flat along $Y$ at $x$.\endproclaim
\demo{Proof} It is well known  that 
the nonsingular points of $Y$
form an open nonempty set
and that  $X$ is normally flat along $Y$ at the points of an open nonempty
subset of $Y$ [7, Corollary (24.5)].\qed \enddemo

\head 2. Ordinary closed points of one dimensional affine reduced schemes.\endhead
In this section $A$ is
 a  reduced local one-dimensional ring of multiplicity $e=e(A)$ and maximal ideal
$\frak m$. With  $K$ we denote the algebraic closure of the residue field
 $k(\frak m)$ of
$A$ and with $\overline A$ we denote the normalization of $A$.\medskip
 
\proclaim{Lemma 2.1} (a) There is a natural immersion $G(A)\subset G(A)\otimes_{k(\frak m)}
 K$; \medskip
 (b) $G(A)\otimes_{k(\frak m)} K$ is  one-dimensional;\medskip
 (c) the Hilbert functions of $A$, $G(A)$ and  $G(A)\otimes_{k(\frak m)} K$
are the same, moreover $e(A)=e(G(A))=e(G(A)\otimes_{k(\frak m)} K)$ and
 $emdim(A)=emdim(G(A))=emdim(G(A)\otimes_{k(\frak m)}  K)$. 

\endproclaim

\demo{Proof}  $(a)$  The immersion follows by the flatness of $G(A)$ 
over $k(\frak m)$\medskip 
$(b)$ [6, $N^0$ 24, Corollaire 4.1.4]\medskip
$(c)$ $G(A)\otimes K$ is the $K$-vector space obtained extending to $K$
the field of scalars of the $k(\frak m)$-vector space $G(A)$, then the Hilbert functions of
$G(A)$ (that is of $A$) and of $G(A)\otimes K$ are the same. This implies the
equalities of the multiplicities and of the embedding dimensions.\qed\enddemo
Since, by Lemma 2.1, $G(A) \otimes_{k(\frak m)}  K$ is a graded one-dimensional algebra
over the algebraically closed field $K$, $Spec((G(A) \otimes_{k(\frak m)}  K)_{red})$ is a 
union of lines each corresponding to a point of $Proj((G(A)\otimes_{k(\frak m)}  K))$.
\definition {Definition 2.2} the lines of $Spec((G(A) \otimes_{k(\frak m)}  K)_{red})$
are the {\sl geometric tangents} of $Spec(A)$ at its closed point.  \enddefinition
Let $n$ be a positive integer and $\frak m^n:\frak m^n=\{b\in \overline A \mid b\frak m^n\subset \frak m^n\}$.
We recall that the subring  of $\overline A$,  $B=\cup_{n>0}(\frak m^n:\frak m^n)$,
is the ring obtained by blowing up the maximal ideal $\frak m$ of A.

 The natural isomorphism
 of schemes $Spec(B/\frak mB) \cong Proj(G(A))$ [6, $N^0$ 32, Lemma 19.4.2] induces an isomorphism of
 schemes\medskip
 $$\phi: Spec(B/\frak mB\otimes_{k(\frak m)}  K)\longrightarrow
Proj(G(A)\otimes_{k(\frak m)}  K)$$ 
By Definition 1.1 a geometric branch $\gamma$ of $Spec(A)$ is
a point of $Spec(\overline A/\frak m \overline A \otimes_{k(\frak m)}  K)$.
The inclusion $B\subset \overline A$ induces a natural homomorphism of zero-dimensional
rings  
$$B/\frak mB\otimes_{k(\frak m)}K \longrightarrow 
\overline A/\frak m \overline A\otimes_{k(\frak m)}K$$
and then a morphism
$$\psi:Spec(\overline A/\frak m \overline A\otimes_{k(\frak m)}K)\longrightarrow
Spec(B/\frak mB\otimes_{k(\frak m)}K)$$ Hence
$\phi(\psi(\gamma))$ is a point of $Proj(G(A)\otimes_{k(\frak m)}  K)$.
\definition{Definition 2.3} The line of $Spec((G(A) \otimes_{k(\frak m)}  K)_{red})$
corresponding to $\phi(\psi(\gamma))$ is
 a geometric tangent of $Spec(A)$ which we call
the {\sl geometric  tangent to $\gamma$}.\enddefinition
\definition{Definition 2.4} The closed point of the scheme $Spec(A)$
is said to be {\sl  ordinary} if $Proj(G(A) \otimes_{k(\frak m)}K)$
is reduced.\enddefinition 

\proclaim {Theorem 2.5} The following conditions are equivalent:\medskip
$(a)$ the closed point of $Spec(A)$ is ordinary\medskip
$(b)$ $Proj(G(A)\otimes_{k(\frak m)}  K)$ consists of $e$ points\medskip
$(c)$ there are $e$ geometric tangents to $Spec(A)$\medskip
$(d)$ $Spec(A)$ has $e$ linear geometric branches with distinct geometric tangents.\endproclaim
\demo{Proof} $(a) \Leftrightarrow  (b)$. Set $R=G(A) \otimes_{k(\frak m)}  K$
and let $\frak q_i$, $i=1,...,n$, be 
the minimal primary ideals of $R$ and $\frak p_i=\sqrt{\frak q_i}$ 
be the corresponding minimal
primes. We have to prove that $Proj(R)$ is reduced if and only if $n=e=e(A)$.
But $Proj(R)$ is reduced if and only if $\frak q_i=\frak p_i$, for any $i$.
By standard facts on multiplicity we have 
$e(R)=\sum_{i=1}^n\lambda(\frak q_i)e(R/\frak p_i)=\sum_{i=1}^n\lambda(\frak q_i)$,
where $\lambda(\frak q_i)$ denotes the lenght of the ideal $\frak q_i$
(note that $e(R/\frak p_i)=1$ because $Spec(R/\frak p_i)$ is a line). Moreover, by Lemma 2.1,(c), 
$e(A)=e(R)$ and then the claim follows from the fact that $\lambda(\frak q_i)=1$
if and only if $\frak q_i$ is prime.\medskip
$(b) \Leftrightarrow (c)$. Follows immediately from Definition 2.2.\medskip
$(c) \Leftrightarrow (d)$. By definition each geometric branch has a geometric tangent
and then the claim is clear if we consider Proposition 1.4,(c)\qed\enddemo 
\example{Example 2.6} Let the normalization $\overline A$ be finite over $A$. If 
 $k(\frak m)$ has characteristic zero and 
$A$ is a seminormal
ring the closed point of $Spec(A)$ is an ordinary point with $e(A)=emdim(A)$ tangents.
 In fact  $A$ is seminormal if and only if $emdim(A)=e(A)$ and  $G(A)$ is reduced [3, Theorem 1]. But,
By Lemma 2.1,(a) we have the inclusion $G(A)\subset G(A)\otimes_{k(\frak m)} K$. Hence  
$G(A)\otimes K$ is the extension to $K$ of the $k(\frak m)$ algebra $G(A)$ and is 
reduced if  $k(\frak m)$ has characteristic zero. In general, if $A$ is seminormal,
 the closed point of $Spec(A)$
needs not to be ordinary. In fact let $k$ be a field of characteristic two with algebraic 
closure $\overline k$. Let $b^2=a=-a\in k$, with $b\in \overline k-k$. Set
$A=(k[X,Y,Z]/(Y^2-aX^2-X^3))_{(X,Y,Z)}$. Clearly $e(A)=emdim(A)=2$ and 
$G(A)=k[X,Y]/(Y^2-aX^2)$ is reduced. Hence $A$ is seminormal. But $Y^2-aX^2=(Y+bX)^2$. Moreover
$k(\frak m)=k$ and $K=\overline k$. Then $G(A) \otimes_{k(\frak m)} K=
\overline k[X,Y]/(Y+bX)^2$ is not reduced.\endexample
Let $R=k[x_0,...,x_r]$ be a reduced standard graded finitely generated $k$-algebra
 over an algebraically closed field $k$. Let $dim(R)=1$ and set $e(R)=e$. Then 
$Proj(R)=\{P_1,...,P_e \}\subset\Bbb P^r$ is a finite set of points . Vice versa
any set of projective points $\{P_1,...,P_e \}\subset\Bbb P^r$ 
is equal to $Proj(R)$ where $R$ is
 its homogeneous coordinate ring. It is well known that, for any $n\in\Bbb N$,
 $H^0(R,n) \leq Min\{e,$$n+r\choose r$$\}$. 

\definition {Definition 2.7} The set  
$\{P_1,...,P_e \}\subset\Bbb P^r$
is {\sl in generic position} in $\Bbb P^r$ (or the points $P_1,...,P_e$  are
  {\sl in generic position} in $\Bbb P^r$)
if the Hilbert function of
$R$ is {\sl maximal} that is $H^0(R,n) = Min\{e,$$n+r\choose r$$\}$  [11, Definition 3.1].\enddefinition
\remark {Remarks} 1. It is proved in [4, Theorem 4] that, for any $e$ and $r$, 
 \lq\lq generic position\rq\rq\ is an open nonempty condition.

2. For details on the notion of
points in generic position.
see  [11]\endremark

\example {Example 2.8} It is easily seen that
any set of  
points of $\Bbb P^1$ is in generic position.\endexample
\example {Example 2.9} A set of  $n+r\choose r$ points in $\Bbb P^r$ $(n>0, r>0)$
is in generic position if and only if they do not lie on a hypersurface of degree
$n$ [11, Corollary 3.4], in particular six points in $\Bbb P^2$ are in generic position
if and only if they do not lie on a conic.\endexample

 Let $S=R[T_1,...,T_n]$, $n\geq 0$, be the polynomial ring over a reduced standard graded
finitely generated $k$-algebra $R=k[x_0,...,x_r]$ ($k$ algebraically closed). Let $dim(R)=1$.
Then
$Spec(S)=\{L_1,...,L_e \}\subset\Bbb A^{n+r+1}$, is  a set of linear varieties 
of dimension $n+1$
all containing a linear subvariety (through the origin)
$L$ of dimension $n$. Vice versa any set of such linear varieties is
equal to $Spec(S)$  where $S$ is its affine coordinate ring.
\definition {Definition 2.10} The set of linear varieties $\{L_1,...,L_e \}$  
is {\sl in generic position} in  $\Bbb A^{n+r+1}$ if  the set of points
 $Proj(R)$ is in generic position in $\Bbb P^r$.\enddefinition

\example {Example 2.11}
 A union of lines of $\Bbb A^{r+1}$ through the origin is 
in generic position if 
the corresponding set of projective 
points of $\Bbb P^r$ is in generic position. In particular, if $emdim(A)=r+1$,
 the set of geometric tangents
$Spec((G(A) \otimes_{k(\frak m)}  K)_{red})\subset \Bbb A^{r+1}$ of  $Spec(A)$ at its maximal ideal is in
generic position in $\Bbb A^{r+1}$ if the Hilbert function of $(G(A) \otimes_{k(\frak m)}  K)_{red}$
is maximal.\endexample
\proclaim{Theorem 2.12} Let $emdim(A)=r+1$. If the closed point of $Spec(A)$ is ordinary with
 geometric tangents
in generic position  in $\Bbb A^{r+1}$ then  the
rings  $G(B)\otimes _{k( \frak m )}  K$ and $G(B)$ are reduced. 
\endproclaim

\demo{Proof} Set $D=G(A)\otimes K$ and $D_{red}=D/nil(D)$. By the surjective graded
homomorphism $\phi:D\rightarrow D_{red}$ we deduce that $H(D_{red},n)\leq H(D,n)$,
for any $n$. If we prove that $H(D_{red},n)= H(D,n)$ then $\phi$ is an isomorphism and
$D$ is reduced (and $G(A)\subset D$ [Lemma 2.1,(a)] is reduced).
 But $H(D_{red},n)\leq H(D,n)\leq Min\{e,$$n+r\choose r$$\}$, and 
$H(D_{red},n)=Min\{e,$$n+r\choose r$$\}$, by assumption, whence the result.
\qed\enddemo
\remark{Remark} If $A$ is the local ring  at a point of a curve over an algebraically
closed field $k$ then $k(\frak m)=k$ and $G(A) \otimes_{k(\frak m)}  K=G(A)$.
Hence, in this case Theorems 2.5 and 2.12 give the results of [10, Lemma-Definition 2.1 and
Theorem 3.3]. Although, by Theorem 2.12, an ordinary point of a curve 
has, in general, reduced tangent cone, there are large classes of examples of ordinary points whose
tangent cone is not reduced [9, Section 4].\endremark
 \head 3. Ordinary codimension one subvarieties.\endhead
In this section
 $X=Spec(R)$ is a reduced variety over an algebraically closed
field $k$ and $Y=Spec(R/\frak q)$ is an irreducible   codimension one subvariety $Y$ 
of $X$ (i.e.  $R_{\frak q}$ is one-dimensional). By $K$ we  
denote 
the algebraic closure of the residue field
 $k(\frak q )$ of $R$ at $\frak q$. We say that  
$e(R_\frak q)=e$ is the {\sl multiplicity} of $Y$ (on $X$) and we set
$emdim(R_\frak q)=r+1$. Note that, if $A$ is the local ring of $X$ at a closed point 
$x$ of $Y$ and $\frak p=\frak qA$
is the prime ideal in $A$ defining $Y$, then $A_\frak p=R_\frak q$. Hence $k(\frak q)$
is equal to the
residue field $k(\frak p)$ of $A$ in $\frak p$.
\definition {Definition 3.1} Let $A$ be the local ring of $X$ at a closed point $x$ of $Y$.
Let $\gamma$ be a  branch of $X$ at  $x$  and $\frak p$ the corresponding
minimal prime of $A^h$ [Lemma 1.2].    $Spec(G(A^h/\frak p))$   is  {\sl the tangent cone
to the  branch $\gamma$} and $Spec(G(A^h/\frak p)_{red})$   is the tangent cone, 
{\sl as a set}, to
$\gamma$. If the branch $\gamma$ is linear the ring $A^h/\frak p$ is regular 
[Definition 1.3] and then its tangent
cone coincides with the tangent space and we call it the {\sl tangent space} to
the branch $\gamma$.\enddefinition
\proclaim {Theorem 3.2} Let $x$ be a point of $Y$ and $A$ be the local ring of $X$ at $x$.
Suppose $Y$ is nonsingular at $x$ and    
$X$ is normally flat along $Y$ at $x$. Let $\gamma$ be a branch of $X$ at $x$. Then
$Spec(G(A)_{red})$ ( that is the tangent cone to $X$ at $x$, as a set)
 is union of linear varieties.
Moreover the tangent cone
to the  branch $\gamma$  is, as a set,  a
 linear variety of $Spec(G(A)_{red})$ and vice versa each linear variety of 
$Spec(G(A)_{red})$ is the tangent cone, as a set,
of at least one branch.\endproclaim
\demo{Proof} $Spec(G(A)_{red})$ is union of linear varieties if and only if the minimal
primes of $G(A)_{red}$ are generated by linear forms.
Moreover, if $\frak p$ is the prime of $Y$ in $A$,
 $(G_\frak p (A) \otimes _{A/ \frak p} k)_{red}$   is  a one-dimensional 
reduced $k$-algebra over the algebraically closed field $k$
 and then its minimal primes are
 generated by linear forms. But the isomorphism
  $G(A) \cong (G_\frak p (A) \otimes _{A/ \frak p} k)[T_1,...,T_d]$
of Theorem 1.6,(a) induces an isomorphism of graded $k$-algebras
 $G(A)_{red} \cong (G_\frak p (A) \otimes _{A/ \frak p} k)_{red}[T_1,...,T_d]$
 and   then also   the minimal primes of  $G(A)_{red}$ are generated by linear forms. 

The second claim
 is proved in [2,Theorem 4.5].\qed\enddemo 

\definition {Definition 3.3} Let $Y$ have multiplicity $e$ on $X$.\medskip

(a) A closed point $x$ of $Y$ is an {\sl ordinary point}
 if, as a set, the tangent cone to $X$ at $x$
is the union of $e$ linear varieties;\medskip
(b) $Y$ is an {\sl ordinary subvariety} of $X$ if the closed point of 
the one-dimensional scheme $Spec(R_\frak q)$ is
ordinary.\enddefinition
   
\proclaim {Theorem 3.4} Let  $Y$ have multiplicity $e$ on $X$ and $x$ be a closed point of $Y$.
If  $Y$ is nonsingular at $x$ and    
$X$ is normally flat along $Y$ at $x$, then $x$ is an ordinary point of $X$ 
if and only if $X$ has, at $x$, $e$ linear branches with distinct tangent spaces.
Moreover if these tangent spaces are in generic position then the tangent cone of $X$ at
$x$ is reduced.\endproclaim
\demo{Proof} Let $A$ be the local ring of $X$ at $x$. 
If  $x$ is an ordinary point of $X$ then $Spec(G(A)_{red}$
is the union of $e$ linear varieties [Definition 3.3] which by Theorem 3.2 are the tangent spaces of $e$
branches and these are linear by   Proposition 1.4,(c). If $x$ has $e$ branches with distinct
tangent spaces then, by Theorem 3.2 these are components of $G(A)_{red}$ and then 
$x$ is ordinary. Now we prove the second claim. Let $\frak p$ be the prime of 
$Y$ in $A$ and consider the isomorphism
$G(A)_{red} \cong (G_\frak p (A) \otimes _{A/ \frak p} k)_{red}[T_1,...,T_d]$
induced by the isomorphism of Theorem 1.6,(a).
By  Definition 2.10  $Spec(G(A)_{red})$
is the union of $e$ linear varieties in generic position if and only if 
$Proj((G_\frak p (A) \otimes _{A/ \frak p} k)_{red})$
 consists of points in generic position that is,
if we set $D=G_\frak p (A) \otimes _{A/ \frak p} k$, $H(D_{red},n)=
Min\{e,$$n+r\choose r$$\}$ But the Hilbert function of 
$G_\frak p (A) \otimes _{A/ \frak p} k$ is equal to the 
Hilbert function of 
$G(A_\frak p)$ [Theorem 1.6,(b)] and then $e=e(A)=e(G(A_\frak p))$. Hence
$H(D,n)\leq Min\{e,$$n+r\choose r$$\}$. Then $H(D_{red},n)=H(D,n)$, 
and $D=D_{red}$ is reduced. Thus  also 
$G(A)$,  which
 is isomorphic to a polynomial ring over $D$, is reduced.\qed\enddemo
\proclaim{Theorem 3.5} The following conditions are equivalent:\medskip
(a) $Y$ is an ordinary subvariety of $X$\medskip
(b) there exists an open nonempty subset $U$ of $Y$ 
such that every closed point $x$ of $U$ is ordinary.\medskip
(c) there exists an open nonempty subset $U$ of $Y$ 
such that $X$ has $e$ linear branches with distinct tangent spaces, at any point of $U$.
\endproclaim 

\demo{Proof} Let $Y$ have multiplicity $e$ on $X$. Then 
$e(G(R_\frak q)\otimes _{k( \frak q )} K)=e$ [Lemma 2.1,(c)].\medskip
$(a)\Rightarrow (b)$
 If $Y$ is ordinary $Proj(G(R_\frak q)\otimes _{k( \frak q )} K)$ has $e$ points
(that is $e$ irreducible components)  
[Definition 3.3,(b) and Theorem 2.5,$(a)\Rightarrow (b)$]. Then by [6, $N^0$ 28, Proposition 9.7.8]
there exists an open nonempty subset $U_1$ of $Y$ such that, for every closed point $x$
of $U_1$, $G_\frak p (A) \otimes _{A/ \frak p} k$
($A$ local ring of $X$ at $x$) has $e$ irreducible components that is $e$ points
[see the proof of Theorem 2.5,$(b)\Leftrightarrow (a)$]. 
Then the $e$ minimal primes of $G_\frak p (A) \otimes _{A/ \frak p} k$
are generated by linear forms. Furthermore, by Theorems 1.6 and 1.8, there exists an open nonempty
subset  $U_2$ of $Y$ such that, for every closed point $x$ of $U_2$,  
 $G(A)$ is isomorphic to a polynomial ring over 
$G_\frak p (A) \otimes _{A/ \frak p} k$. Then if $x$ is a point of $U_1\cap U_2$
 the minimal primes of $G(A)$ are extensions of
 the minimal primes of the one dimensional finitely generated graded $k$-algebra
 $G_\frak p (A) \otimes _{A/ \frak p} k$, hence they are generated by linear forms.\medskip 

$(b)\Rightarrow (a)$ Suppose that, at every point $x$ of an open nonempty subset of $Y$, 
$Spec(G(A)_{red})$ has $e$ irreducible
components, that is $G(A)_{red}$ has $e$ minimal primes. Then, by the graded isomorphism
$G(A)_{red} \cong (G_\frak p (A) \otimes _{A/ \frak p} k)_{red}[T_1,...,T_d]$,
 on the points of an open set of $Y$, $(G_\frak p (A) \otimes _{A/ \frak p} k)_{red}$ has $e$
minimal primes [6, $N^0$ 28, Proposition 9.7.8]. Hence the zero-dimensional scheme 
$Proj(G(A_\frak p)\otimes _{k( \frak p )} K)=
Proj(G(R_\frak q)\otimes _{k( \frak q )} K)$
has  $e$ points, and then, since by Lemma 2.1,(c) $e(G(A_\frak p)\otimes _{k( \frak p )} K)=
e(G(A_\frak p))=e(A_\frak p)=e$, is reduced [Theorem 2.5, $(b)\Rightarrow(a)$].\medskip
$(b) \Leftrightarrow (c)$ It is an easy consequence of Theorem 1.8 and Theorem 3.4.\qed\enddemo

Let $Y$ be an ordinary  subvariety of codimension one and of multiplicity $e$
on a reduced variety $X$. Let $x$ be a closed point of $Y$ such that $Y$ 
is nonsingular at $x$ and $X$ is normally flat along $Y$ at $x$.
Then by Theorem 3.2 the tangent cone of $X$ at $x$ consists, as a set, of  linear
varieties, but the number of these can be less then $e$ and then $x$ is not an
ordinary point of $X$ as the following example shows. 
\example{Example 3.6} Let $R=\Bbb C[X_1,X_2,X_3]/(X_1X_2^n-X_3^n)=
\Bbb C[x_1,x_2,x_3]$ ($n\geq 2$)  and $A$ be the local ring of $X=Spec(R)$ at the 
maximal ideal $(x_1,x_2,x_3)$. 
The non-normal locus of the hypersurface $X$ is the line 
$Y:x_2=0,x_3=0$ of multiplicity $n$ on $X$ and if
$a\in\Bbb C$, $a\neq 0$, the tangent cone
 at the  point $(a,0,0)$ of $L$  
consists, as a set, of the $n$ distinct
planes $x_3=bx_2$ , where $b^n=a$. Then $Y$ is an ordinary subvariety of $X$. Moreover
$e(A)=e(R_\frak q)=n$, where $\frak q=(x_2,x_3)$.
Then, by [7, Corollary (23.22)], $X$ is normally flat along $Y$ at $x$ and clearly $Y$
is non singular at $x$. But the tangent cone of $X$ at $(0,0,0)$
is $Spec(\Bbb C[X_1,X_2,X_3]/(X_3^n)$ and then, as a set, consists of the plane $x_3=0$.
But $e=e(A)=n\geq 2$ and $(0,0,0)$ is not ordinary.\endexample

\enddocument